\documentstyle[12pt]{article}
\begin{document}

\title{Mathematical Support to Braneworld Theory}

\author{Edmundo M. Monte\thanks{E-Mail: edmundo@fisica.ufpb.br and edmundomonte@pq.cnpq.br}\\
Departamento de Fisica-Matematica, \\
Universidade Federal da Paraiba, 58059-970,
\\ Jo\~{a}o Pessoa, Paraiba, Brazil.} 

\maketitle


\begin{abstract}
The braneworld theory appear with the purpose of solving the problem of the hierarchy of the fundamental interactions. The perspectives of the theory emerge as a new physics, for example, deviation of the law of Newton's gravity. One of the principles of the theory is to suppose that the braneworld is local submanifold in a space of high dimension, the bulk, solution of Einstein's equations in high dimension. In this paper we approach the mathematical consistency of this theory with a new proof of the fundamental theorem of submanifolds for case of semi-Riemannian manifolds. This theorem consist an essential mathematical support for this new theory. We find the integrability conditions for the existence of space-time submanifolds in a pseudo-Euclidean space.
\end{abstract}

\section{Brief Historic of Immersions of the Space-times}

One of the main problems of Riemannian geometry is the
immersion of a Riemannian manifold in Euclidean space. This
problem consists to know when a Riemannian manifold $M$, with
metric $g$, permits an immersion $f:M\longrightarrow R^{n}$ such
that the induced metric has its origin in the scalar product of
$R^{n}$ in the submanifold $f(M)\subset R^{n}$ and coincides with
$g$. When $f$ exists we say that $M$ is isometrically immersed
into $R^{n}$.

This problem can be approached under some aspects. We can require
the isometric immersion $f$ to be defined on $M$ (global case) or
only in the neighborhood of some point of $M$(local case). In
general the dimension of the Euclidean space in which a immersion
of $M$ is realized is high, which depends of the regularity of the
immersion. When is possible to release an isometric immersion for
sufficiently large $n$ it another question appears: the least
possible value of $n$. Some other questions appear from immersion
problem, for instance: isometric deformation, extrinsic absolute
invariants, algebraic criteria for immersion class, rigid
immersion, etc. Nowadays there is a great literature about
immersion, embedding (imbedding), and submanifolds connected with
others branches of mathematics.

The oldest publication cited in contemporaneous articles about
isometric immersion is of 1873 by Schlaefli, where he made a
conjecture: Any Riemannian manifold $m$ dimensional can be
isometrically immersed into Euclidean space $D=\frac{m(m+1)}{2}$
dimensional. In 1926, Janet published a proof of the conjecture,
but in 1931 Burstin gave a correct proof. Cartan gave a correct
proof too, but used differential forms. There are many results in
the literature, some of them are: Nash's theorem (global
immersions), O'Neill (algebraic criteria for immersion),
Friedman's theorem (generalization of Bustin's theorem for
pseudo-Riemannian manifolds and analytical metrics), Greene's
theorem (theorem for pseudo-Riemannian manifolds and differential
metrics), etc.

The physics of gravitational fields is intimately connected with
the geometry of four-dimensional pseudo-Riemannian manifold
(space-times) whose metrics obey Einstein's field equations. The
Einstein's theory customarily has been treated as an intrinsic
geometric theory where the four-dimensional space-time consists of
an arena of the physics events perfectly according with
pseudo-Riemannian geometry. We live in a four dimensional world.
However the idea of unification of fundamental interactions leads
us to higher dimensional theories, where now the arena of physical
events passes to be a high dimension space and without doubt we
must have our world immersed into the latter.

The immersion problem which involves pseudo-Riemannian manifolds
is brought to multidimensional physics where it may help in
finding new geometrical characterizations of gravitational fields
which can be connected to physics. In submanifolds (space-times)
this problem is called embedding of space-times.

The first result about embedding of space-times came from the
paper of Kasner, 1921, "Finite Representation of the Solar
Gravitational Field in Flat Space of Six Dimensions". He gave an
explicit local isometric embedding of the external Schwarzschild
solution. Others results about this topic have been found in the
literature with various physical motivations, or simply for
mathematics motivations.

In 1926 Kaluza and Klein tried to explain both gravitational
interaction and the electromagnetic interaction in one theory.
They worked in a five dimensional space-time, where one of the space
dimensions is compact. It was extended to include all the gauge
interactions. Thus to accommodate the strong, weak and
electromagnetic interactions we need at least seven extra
dimensions plus four dimensions of our world. Only in the eighties
Kaluza-Klein theory was derived from the hypothesis that the
four-dimensional space-time is locally and isometrically embedded
in higher dimensional space.

In 1965 occurred a great event with mathematicians and physicists
interested on the embedding problem motivated in properties of
particle physics in curved space-time, identifying the internal
symmetries of particle physics with the symmetries of the normal
piece of embedding space.

The embedding problem has been required in problems linked to
minimal class of the embedding, extrinsic gravity, theory of
strings and membranes. Nowadays the embedding problem emerged with
a new theory - braneworld - a theory that lead us the idea of
unification of fundamental interactions using extra dimensions.
Such a model has been showing positive in the sense that we find
perspectives and probably deep modifications in the physics, such
as: unification in a TeV scale, quantum gravity in this scale and
deviation of Newton's law of gravity for small distances. 

A braneworld may be regarded as a space-time locally embedded in a
higher dimensional space, the bulk,  solution of higher
dimensional Einstein's equations. Furthermore, the embedded
geometry  is assumed to  exhibit quantum  fluctuations with
respect to  the extra dimensions at the  TeV scale of  energies.
Finally, all  gauge interactions belonging to the standard model
must remain confined to the four-dimensional space-time.
Contrasting with  other  higher dimensional theories, the extra
dimensions may be large and even infinite, with the possibility of
being observed by  TeV accelerators. The embedding  conditions
relate the  bulk  geometry  to the braneworld geometry, as it is
clear from the  Gauss-Codazzi-Ricci equations. \cite{Nash}

The fundamental theorem of Riemannian geometry warrants that a given 
Riemannian manifold can be locally and isometrically immersed in 
higher dimensional Euclidean spaces, if and only if there exist tensors that 
satisfy the Gauss-Codazzi-Ricci equations. \cite{Jacob,Keti,Einsenhart}

However we give a new proof of a fundamental and classical theorem for manifolds 
with  arbitrary signatures and dimensions. This theorem is very interesting when 
we consider a space-time of the General Relativity as a submanifold of a pseudo-Euclidean space.  

Braneworlds  may be  described  as  a family of stable  perturbations  of  a  given locally embedded  background  space-time. This suggests having a dynamic bulk whose geometry depends on that of the braneworld. From the point of view of perturbations in braneworlds we can suppose that there is a fundamental state (ground state) of the bulks named a induced pseudo-Euclidean bulk so that any braneworld would be immersed in some appropriate flat space bulk. This immersion had been guaranteed by the fundamental theorem of submanifolds.\cite{Edmundo,Edmundo1}

\section{Preliminaries}

Consider a point $p \in (V_n,g)$, where $(V_n,g)$ is a pseudo Riemannian 
manifold (braneworld) with $r-s$ signature, given by local coordinates $x^i, i = 1, \ldots, n$. 
We have that the local and isometric immersion $Y : (V_n, g) \to (V_D, \eta)$, where 
$(V_D,\eta)$ is a pseudo Euclidean manifold (bulk) with $p-q$ signature, the 
position of $p \in V_n$ is indicated for $Y^\mu =Y^\mu(x^1, \ldots, x^n)$. In 
local coordinates the induced metric is given by 
\begin{equation}
g_{ij} = \eta_{\mu\nu}Y^\mu_{, i} Y^\nu_{, j},
\end{equation}
where $Y^\mu_{, i} = \frac{\partial Y^\mu}{\partial x^i}$; the local orthonormal 
basis of $T_pV_n$ is $\{\frac{\partial}{\partial x^i}\}$ and another local 
orthonormal basis of $T_{Y(p)}V_D$ is $\{\frac{\partial}{\partial Y^\mu}\}$.

Set $\nabla'$ the Levi-Civita connection of $V_D$, where  $U'$ and $W'
\in V_D$. We decompose $\nabla'_{U'}W'$ on $V_D$, when we restrict the 
differential fields $U'|_U$ and $W'|_W$, on the fiber $TV_n$. The  
decomposition of $\nabla'_{U}W$ is given by
\begin{equation}
\nabla'_{U}W = \nabla_UW + b(U,W),
\end{equation}
where the tangent part $\nabla_UW$ is the projection of $\nabla'_UW$ on the 
$TV_n$. The normal part $b(U,W)$ is defined by eq. (2), it is a 2-linear symmetric 
application. 
\[
b: TV_n \times TV_n \to T(V_n)^\perp,
\]
we called it the second fundamental form of $V_n$. Observe that $\nabla_UW$ 
is a Levi-Civita connection of $V_n$,  agreeing with the Levi-Civita 
connection of the induced metric. We have that $T(V_n)^\perp$ is the orthogonal 
complement of $TV_n$ in $TV_D | _{V_n}$. \cite{Kobayashi}

There exists also an induced connection on the fiber $T(V_n)^\perp$. This 
connection is defined by the application, 
\[
\nabla'':TV_n \times T(V_n)^\perp \to T(V_n)^\perp
\]
with
\begin{equation}
\nabla'_UN = \nabla''_UN+ S(U, N)
\end{equation}
the decomposition of $\nabla'_UN$, where the tangent part is given by the 
linear application, 
\[
S:TV_n \times T(V_n)^\perp \to TV_{n}
\]
related with $b$ through
\begin{equation}
\eta(S(U,N),N) = \eta(b(U,W),N)
\end{equation}
while $\nabla''_UN$ is the normal part of $\nabla'_UN$ in eq. (3). \cite{Kobayashi}

In local coordinates we can write, 
\[
U = u^i \frac{\partial}{\partial x^i}, \quad W = w^j
\frac{\partial}{\partial x^j}, \quad N_A =
N^\mu_A\left(\frac{\partial}{\partial Y^\mu}\right)^\perp,
\]
where $\{\frac{\partial}{\partial Y^\mu}\}^\perp , \; \mu = D-n,
\ldots, D$ is local and orthonormal basis of $T_p(V_n)^\perp $, with 
\begin{equation}
\eta(N_A, N_B) = \eta_{AB} = \epsilon_A g_{AB}\;\;(\mbox{com}\;\; \epsilon_A
= \pm 1)
\end{equation}
and
\begin{equation}
\eta(N_A, U) = 0.
\end{equation}
We indicate $N_A$ any normal vector. We have that,
\begin{equation}
\nabla_UW = u^i(w^k_{, i} + \Gamma^k_{ij}w^j)
\frac{\partial}{\partial x^k},
\end{equation}
\begin{equation}
\nabla''_UN_A = u^i(N_{A,i}^\mu + A^\mu_{i\nu}
N_A^\nu)(\frac{\partial}{\partial Y^\mu})^\perp .
\end{equation}
where $\mu,\;\;A, B, \ldots = D - n, \ldots, D.$ Thus,  
$A_{i\mu\nu} \equiv A_{iAB}$ and $b\left(\frac{\partial}{\partial x^i}, \frac{\partial}{\partial x^j} \right) = b_{ij}
\mbox{\phantom{)}}^A(\frac{\partial}{\partial Y^A})^\perp$. The second form 
coefficients $b_{ijA}$ are the components of $D-n$ tensors on $T_pV_n$, for 
each fixed $A$, with $b_{ijA} = b_{jiA}$. On the other hand, $A_{iAB}  =
-A_{iBA}$ are components of 
$\left(\begin{array}{c} D - n\\ 2 \end{array} \right)$ vectors on $T_pV_n$, it corresponds to the    
normal part of the connection on $V_D$ and we called it twisting vector.\cite{Gonner}

Replacing $U = \frac{\partial}{\partial x^i}$, and $\; N_A
= \left(\frac{\partial}{\partial Y^A}\right)^\perp$ in the equations (2) and (3) 
we obtain,
\begin{equation}
Y^\mu_{; ij} = g^{AB}b_{ijA}N^\mu_B
\end{equation}
\begin{equation}
N^\mu_{A, j} = -g^{mn}b_{jmA}Y^\mu_{, n} + g^{MN}A_{jAM}N^\mu_N.
\end{equation}
where $Y^\mu_{; i}$ is the covariant derivate of $Y^\mu$ with respect the 
$g_{ij}$ metric. The equations (9) and (10) are called: Gauss formula and 
Weingarten formula, respectively.

It is easy to see from eq. (9) that the second form coefficients are given by 
\begin{equation}
b_{ijA} = Y^\mu_{; ij} N^\nu_A \eta_{\mu\nu}.
\end{equation}
>From eq. (10) we obtain the twisting vector components
\begin{equation}
A_{iAB} = N^\mu_AN^\nu_{B, i} \eta_{\mu\nu}.
\end{equation}
We have that the system of partial differential equations from eqs. (9) and (10) 
can be resolved for arbitrary initial vectors of $\{Y^\mu_{, i} , N^\mu_A \}$ 
for a fixed point, so that they satisfy,
\begin{equation}
g_{ij} = Y^\mu_{, i} Y^\nu_{, j} \eta_{\mu\nu},
\end{equation}
\begin{equation}
N^\mu_AY^\nu_{, i} \eta_{\mu\nu} = 0,
\end{equation}
\begin{equation}
N^\mu_AN^\nu_B\eta_{\mu\nu} = g_{AB} = \epsilon_A\delta_{AB},
\end{equation}
where $\epsilon_A = \pm 1$ and depends from the $(V_D,\eta )$ signature.

In the next section we are going to prove the fundamental theorem of the submanifolds. Particularly in the case of braneworlds (space-times) submanifolds the theorem appear as mathematical support to braneworld theory.

\section{The Fundamental Theorem of the Submanifolds}


{\it A pseudo Riemannian manifold $(V_n,g)$ with $r-s$ signature, is locally and 
isometrically immersed in the pseudo Euclidean manifold $(V_D,\eta )$ with $p-q$  
signature, if only if there exist $D-n$ symmetric matrices $b_{ijA}$ and 
$\left(\begin{array}{c} D - n\\ 2 \end{array} \right)$ vectors $A_{iAB}$ satisfying the Gauss-Codazzi-Ricci 
equations.}

\vspace{3mm}
{\bf Proof:}
\vspace{3mm}

Suppose that there exists a local and isometric immersion 
\[
\begin{array}{lccl}
Y: & (V_n, g) & \to & (V_D, \eta)\\
\vspace{4mm} 
& x^i & \to & Y^\mu = Y^\mu(x^i)
\end{array}
\]
If to a fixed point we choose arbitrary tangent vectors 
$Y^\mu_{,i}$ and normal vectors $N^\mu_{A}$, then these vectors 
satisfy eqs. (13), (14) and (15) from section (2). On the other hand, the set 
$\{Y^\mu_{, i} ,N^\mu_{A}\}$ needs to satisfy the Gauss and Weingarten formulas 
given by the system 
\begin{equation}
Y^\mu_{,ij} = g^{AB}b_{ijA}N^\mu_{B} - \Gamma^r_{ij}Y^\mu_{,r}
\end{equation}
\begin{equation}
N^\mu_{A,j} = -g^{ml}b_{jmA}Y^\mu_{,l} + g^{MN}A_{jAM}N^\mu_{N},
\end{equation}
where $b_{ijA} = Y^\mu_{,ij}N^{\nu}_{A}\eta_{\mu\nu}$
and $A_{iAB} = N^\mu_AN^\nu_{B, i}\eta_{\mu\nu},$ warranted the existence 
of these tensors. This system is super determined and since $Y^{\mu}$ 
already exists the system (1) and (2) has solution.
Therefore the integrability conditions
\begin{equation}                                           
Y^\mu_{i,jk} = Y^\mu_{i,kj}
\end{equation} 
\begin{equation}
N^\mu_{A,jk} = N^\mu_{A,kj}
\end{equation}
are satisfied for $D-n$ matrices $(b_{ijA})$ and $\left(\begin{array}{c} D - n\\ 2 \end{array} \right)\;\;$ 
$A_{iAB}$ vectors. These integrability conditions, consequence of Frobenius 
theorem, conduce to the Gauss-Codazzi-Ricci equations.\cite{Boothby}

Using (1) calculate (3) and use the Weingarten formula for 
$N^\mu_{B, k}$ and $N^\mu_{B, j}$. We obtain,

\[
\left\{\begin{array}{ll}
g^{AB}b_{ijA, k}N^\mu_B + &
 g^{AB}b_{ijA} \left(-g^{m\ell}b_{kmB}Y^\mu_{,\ell} +
g^{MN}A_{kBM}N^\mu_N\right)  
\vspace{2mm}\\
&- \Gamma^r_{ij, k}Y^\mu_{, r} - \Gamma^r_{ij}Y^\mu_{,rk} +
g^{AB}_{,k}b_{ijA}N^\mu_{B}  = 
\end{array}\right.
\]
\begin{equation}
\left\{\begin{array}{ll}
g^{AB}b_{ikA,j}N^\mu_B +&
 g^{AB}b_{ikA}\left(-g^{m\ell}b_{jmB}Y^\mu_{,l} +
g^{MN}A_{jBM}N^\mu_{N}\right)
\vspace{2mm}\\
& - \Gamma^r_{ik} Y^\mu_{,rj} +
g^{AB}_{,j}b_{ikA}N^\mu_B 
\end{array}\right.
\end{equation}

Using that $\{x^i\}$ are geodesic coordinates, then $\Gamma^k_{ij} = 0$. 
Multiply (5) by $\eta_{\mu\nu}Y^\nu_{,n}$ and employ the eqs. (13), (14) and (15) from section (2) to 
obtain the Gauss equation, 
\begin{equation}
R_{nijk} = g^{AB}(b_{ikA}b_{jnB} - b_{ijA}b_{knB})
\end{equation}  
 
Again using that  $\{x^i\}$ are geodesic coordinates, multiplying (5) by 
$\eta_{\mu\nu}N^\nu_C$ and employing again eqs. (13), (14) and (15) from section (2) we obtain, 
\begin{equation}
b_{ijC,k} - b_{ikC,j}  = g^{AB}(b_{ikA}A_{jBC} - b_{ijA}A_{kBC}).
\end{equation}
We have a geodesic system, then we can write the Codazzi equation,
\begin{equation}
b_{ijC;k} - b_{ikC;j} = g^{AB}(b_{ikA}A_{jBC} - b_{ijA}A_{kBC}).
\end{equation}
Now take (2), calculate (4) and use the Weingarten formula for 
$N^\mu_{N,k}$ and $N^\mu_{N,j}$. We obtain,
\begin{equation}
\begin{array}{ll}
- &  g^{m\ell ,k}b_{jmA}Y^\mu_{,\ell} - g^{m\ell}b_{jmA,k}Y^\mu_{,\ell}
- g^{m\ell}b_{jmA}Y^\mu_{,\ell k} + g^{MN}_{,k}A_{jAM}N^\mu_N
\vspace{2mm}\\
+ & g^{MN}A_{jAM,k}N^\mu_N +
g^{MN}A_{jAM}\left(-g^{m\ell}b_{kmN}Y^\mu_{,\ell} +
g^{CD}A_{kNC}N^\nu_{D}\right) = \vspace{2mm}\\
- & g^{m\ell}_{,j}b_{kmA}Y^{\mu}_{,\ell} -
g^{m\ell}b_{kmA,j}Y^\mu_{,\ell} - g^{m\ell}b_{kmA}Y^\mu_{,\ell j}
- g^{MN}_{j}A_{kAM}N^\mu_{N} \vspace{2mm}\\
+ & g^{MN}A_{kAM,j}  N^\mu_{N} +
g^{MN}A_{kAM} \left(-g^{m\ell} b_{jmN}Y^\mu_{,\ell} +
g^{CD}A_{jNC}N^\mu_D\right)
\end{array}
\end{equation}
multiplying (9) by $\eta_{\mu\nu}N^\nu_B$, employing again eqs. (13), (14) and (15) from section (2) and using the definition of the second form we obtain the Ricci equation,
\[
A_{jAB;k} - A_{kAB;j}= 
\]
\begin{equation}
g^{MN}\left(A_{jAM}A_{kNB} -
A_{kAM}A_{jNB}\right) + g^{ml}\left(b_{kmA}b_{ljB} - b_{jmA}b_{lkB}\right) 
\end{equation}

We conclude that if there exists an immersion, then there exist $(D-n)$ 
symmetric $b_{ijA}$ and $\left(\begin{array}{c} D - n\\ 2 \end{array} \right)\;\;$ $A_{iAB}$ vectors 
satisfying the Gauss-Codazzi-Ricci equations.

Now we prove the sufficiency condition. We know that if for a fixed point we  
choose arbitrary vectors $Y^\mu_{,i}$ and $N^\mu_A$ satisfying eqs (13), (14) and (15) section (2), then $Y^{\mu}$ is determined. Therefore we need to calculate 
$\{Y^\mu_{,i}N^\mu_A\}$.  For this realization consider the system equations 
of (1) and (2), where we know to exist by hypothesis $(D-n)$ 
$(b_{ijA})$ and $(1/2)(D-n)(D-n-1)\;\;$ $A_{iAB}$. We observe that the 
system (1) and (2) is a system of partial differential equations super determined. 
The Frobenius theorem warrants the integrability of (1) and (2), for this we 
need to prove that
\begin{equation}
Y^\mu_{i,jk} = Y^\mu_{i,kj}
\end{equation}
and
\begin{equation}
N^\mu_{A,jk} = N^\mu_{A,kj}
\end{equation}
On the other hand, if (11) and (12) are satisfied, the system (1) and (2) 
possesses solution $\{Y^\mu_{,i}, N^\mu_{A}\}$.
By hypothesis we know that the Gauss-Codazzi-Ricci equations are satisfied 
for $b_{ijA}$ and $A_{iAB}$ values. Therefore:\\

- Take (6), using the fact that we can have $\Gamma' s = 0$ on the some 
point, we rewrite (6). Then we have for all $Y^\nu_{,n}$ the projections on 
$T_{p}V_n$ of the vectors $Y^\mu_{,ijk}$ and $Y^\mu_{,ikj} \in T_{p}V_D$ are 
the same, 
\begin{equation}
Y^\mu_{,ijk}Y^\nu_{,n}\eta_{\mu\nu} = Y^\mu_{,ikj}Y^\nu_{,n}\eta_{\mu\nu},
\end{equation}
where the expressions for $Y^\mu_{,ijk}$ and $Y^\mu_{,ikj}$ are given by (5).

- Now take (8) and using eqs (13), (14) and (15) from section (2) we have for all $N^\mu_B$ that the 
projections on $T_{p}(V_n)^\perp$ of the vectors 
$Y^\mu_{,ijk}$ and $Y^\mu_{,ikj} \in  T_{p}V_D$ are equal, 
\begin{equation}
Y^\mu_{,ijk}N^\nu_{C}\eta_{\mu\nu} = Y^\mu_{,ikj}N^\nu_{C}\eta_{\mu\nu}.
\end{equation}
Wherein $TV_D = TV_n \oplus T(V_n)^\perp$ and the projections vectors 
$Y^\mu_{,ijk}$ and $Y^\mu_{,ikj}$ in the spaces $T_{p}V_n$ and 
$T_{p}(V_n)^\perp$ are equal, we obtain (11).

Now we prove (12). Take (10) and using eqs. (13), (14) and (15) from section (2), on the geodesic system we have for all $N^\nu_{B}$, the projections on $T_{p}(V_n)^\perp$ of 
vectors $N^\mu_{A,jk}$ and $N^\mu_{A,kj} \in T_{p}V_D$  are equal,
\begin{equation}
N^\mu_{A,jk}N^\nu_B\eta_{\mu\nu} = N^\mu_{A,kj}N^\nu_B\eta_{\mu\nu},
\end{equation}  
where the expressions for $N^\mu_{A,jk}$ and $N^\mu_{A,kj}$ are given by (9).

Observe that the projections of these vectors are equal zero,
\begin{equation}
N^\mu_{A,jk}Y^\nu_{,n}\eta_{\mu\nu} = 0 , \; \forall \;Y^\nu_{,n}
\in T_{p}V_n
\end{equation}
\begin{equation}
N^\mu_{A,kj}Y^\nu_{,n}\eta_{\mu\nu} = 0 , \; \forall \;Y^\nu_{,n}
\in T_{p}V_n
\end{equation}
therefore
\begin{equation}
N^\mu_{A,jk}Y^\nu_{,n}\eta_{\mu\nu} = N^\mu_{A,kj}
Y^\nu_{,n}\eta_{\mu\nu}.
\end{equation}
Wherein $TV_D = TV_n \oplus T(V_n)^\perp$ and the projections of vectors 
$N^\mu_{A,jk}$ and $N^\mu_{A,kj}$ in the spaces $T_{p}V_n$ and  
$T_{p}(V_n)^\perp$ are equal, we obtain (12).
\rule{2mm}{2mm}

\section{Conclusions}

In this paper we have especially investigated a mathematical support for braneworld theory with the new proof of fundamental theorem of submanifolds. Through of a mechanism of the successive immersions we always can found the ground state of the bulks, i.e., a pseudo-Euclidean bulk. We can always find a space (bulk or braneworld) which contains locally immersed another bulk.Then the theorem appears to warrant that each braneworld can be immersed in the pseudo-Euclidean bulk. Certainly the geometric entities that appear of that immersion formalism, for example, the second fundamental form and the coefficients of normal connection should be interpreted physically in agreement with the theory of braneworlds. We are preparing an article explaining physically this mechanism.

\section{Acknowledgments} 
I would like to thank Professor Marcos Maia for useful discussions, and to FAPES-ES-CNPq-PRONEX-Brasil.

\end{document}